\documentclass{article}       
\usepackage{amsrefs}
\newtheorem{theorem}{Theorem}

\newtheorem{proposition}{Proposition}
\newtheorem{definition}{Definition}
\newtheorem{remark}{Remark}

\def\ttt{\,\dots\,}

\def\ds{\displaystyle}
\def\a{\alpha}
\def\b{\beta}
\def\g{\gamma}

\def\la{\lambda}
\def\o{\omega}
\def\s{\sigma}
\def\an{\hbox{$A_n$}}
\def\eqref#1{\hbox{$(\ref{#1})$}}

\def\G{\Gamma}
\def\l{\lambda}

\def\be#1{\begin{equation}\label{#1}}
\def\ee{\end{equation}}

\def\email{{\sl e-mail: }}
\date{}

\begin{document}

\title{On the existence of pre-semigeodesic coordinates
}


\author{Irena Hinterleitner,            Josef Mike\v s}

\maketitle

\begin{abstract}
In the present paper we consider the problem of the existence of pre-semigeodesic coordinates on manifolds with affine connection. We proved that pre-semigeodesic coordinates exist in the case when the components of the affine connection are  twice differentiable functions. 
\smallskip

\footnote{The paper was supported by the project IGA PrF 2014016 Palacky University Olomouc and No. LO1408 "AdMaS UP - Advanced Materials, Structures and Technologies", supported by Ministry of Education, Youth and Sports under the „National Sustainability Programme I" of the Brno University of Technology.}
\footnote{I. Hinterleitner, 
             	Dept. of Mathematics, Faculty of Civil Engineering,
							Brno University of Technology, 
							\v Zi\v zkova 17, 602 00 Brno, Czech Republic, 
							\email{hinterleitner.i@fce.vutbr.cz}}          
\footnote{J. Mike\v s, 
              Dept. of Algebra and Geometry, Palacky University, 
              17.~listopadu 12, 77146 Olomouc, Czech Republic, 
              \email{josef.mikes@upol.cz}}
			   
\end{abstract}

\section{Introduction}
\label{intro}
Geodesics are fundamental objects of differential geometry, analogous to straight lines in Euclidean space. A geodesic
is a curve whose tangent vectors in all of its points are parallel.
	Some properties of geodesic lines in  mechanics:
a point mass  without  external influences moves on a geodesic line, another example of geodesics is an ideal elastic ribbon without friction between two points on a curved surface \cite{lagg,le}. Geodesics are of particular importance in general relativity. Timelike geodesics in general relativity describe the motion of inertial test particles.

Let $A_n=(M,\nabla)$ be an $n$-dimensional manifold $M$ with affine connection $\nabla$.
A curve $\ell$ in \an\  is a {\it geodesic} when its tangent vector field remains in the tangent distribution of $\ell$ during parallel transport along the curve or, equivalently
	 if and only if the covariant derivative of its tangent vector, i.e. $\la(t)=\dot{\ell}(t) $ is proportional to the tangent vector
	$
	\nabla_{\la}\la=\rho(t)\,\la,
	$
	where $\varrho$ is some function of the parameter~$t$ of the curve $\ell$.
	
	When the parameter~$t$ of the geodesic is chosen so that $\varrho(t)\equiv 0$, then this parameter is called {\it natural} or {\it affine}. 
A natural parameter is usually denoted by $\tau$.
	\smallskip

With geodesics some  special coordinates are closely associated: 

\centerline{geodesic, semigeodesic and pre-semigeodesic coordinates.}
\smallskip

{\it Geodesic coordinates at a point $p$} and {\it along a curve $\ell$} ({\it Fermi coordinates}) are characterized by  vanishing Christoffel symbols (or components of the affine connection) at  the point $p$ and along the curve $\ell$, respectively.\smallskip

Let us consider a non-isotropic coordinate hypersurface 
${\Sigma}$:\ $x^1=c$ in \mbox{(pseudo-)} Riemannian space $V_n$.
Let us fix some point $(c,x^2,\dots ,x^n)$ on~${\Sigma}$ and construct the geodesic $\gamma$ passing through the point and tangent 
to the unit normal of ${\Sigma}$; $\gamma$ is an $x^1$-curve, it is parametrized by $\gamma(x^1)=(x^1+c, x^2,\dots ,x^n)$ and 
$x^1$ is the arc length on the geodesic. 
Coordinates introduced in this way are called 
{\it semigeodesic coordinates}
in~$V_n$.

It is well known that  the metric  of $V_n$ in semigeodesic coordinates has  the following form:
$
ds^2=e\,(dx^1)^2 + g_{ab}(x)\,dx^a\,dx^b$, $a,b>1, \ e=\pm1$. 
On the other hand this coordinate form of the metric is a sufficient condition for the coordinate system to be semigeodesic. In this case for the Christoffel symbols of the second type follows 
$\G^h_{11}=0$, $h=1,\dots, n$.

Advantages of such coordinates are known since C.F.~Gauss ({\it Geo\-d\" a\-tische Parallelkoordinaten}, \cite[p.~201]{krey}),

 \textit{Geodesic polar coordinates}:
can be also interpreted as a ``limit case'' of semigeodesic coordinates:
all geodesic coordinate lines $\varphi=x^2=\mbox{\rm const.}$ pass through one point called the pole, corresponding to
$r=x^1=0$, 
and lines $r=x^1=\mbox{\rm const}$ are geodesic circles
({\it Geod\" atische Polarkoordinaten}, \cite[pp.~197-204]{krey}).  
\smallskip

Let $A_n=(M,\nabla)$  be  an \hbox{$n$-di}\-men\-sio\-nal manifolds $M$  with the affine connection $\nabla$, dimension $n\geq2$, and let
$U\subset M$ be a coordinate neighbourhood at the point $x_0\in U$.
A couple $(U,x)$ is a coordinate map on \an.

Semigeodesic coordinate systems on surfaces and (pseudo-) Riemannian manifolds are generalized in the following way (Mike\v s, Van\v zurov\'a, Hinterleitner \cite[p.~43]{mvh}):
\begin{definition}\label{def1}\rm
Coordinates $(U,x)$ in \an\ are called  \textit{pre-semigeodesic coordinates} if one system  of~coordinate lines is geodesic and the coordinate is just the natural parameter.
\end{definition}

In a paper by J.~Mike\v s and A.~Van\v zurov\'a \cite{miva} these coordinates were called {\it general Fermi coordinates}, and  the reconstruction of components of the affine connection in these coordinates is shown, if we known a certain number of components of the curvature tensor.

In \cite[p.~43]{mvh}, \cite{miva}  the following theorems were proved.

\begin{theorem}\label{th1}
The conditions \
\ $\G^h_{11}(x)=0$, \ \   $h=1,\dots,n$,
are satisfied in  $(U,x)$ if and only if $(U,x)$ is pre-semigeodesic.
\end{theorem}

\begin{theorem}\label{th2}
The conditions \ $\G^h_{11}(x)=0$, \ \   $h=1,\dots,n$, 
are satisfied in a coordinate map $(U,x)$ if and only if the parametrized curves 
$$\ell{:}\ I\to U, \ \ell(\tau)=(\tau,a_2,\dots,a_n), \tau\in I,\
a_i\in R,\ i=2,\dots,n,
$$ 
are canonically parametrized geodesics of $\nabla_{|U}$,
I is some interval, $a_k$ are suitable constants chosen so that  $\ell(I)\subset U$,
$\G^h_{ij}$ are components of the affine connection $\nabla$, the subset $U\subset M$ is a coordinate neighbourhood of $A_n=(M,\nabla)$.
\end{theorem}

We thought that the existence of this chart is trivial. This problem is obviously more difficult than we supposed.
 This was observed in \cite{dusko1,dusko2} where precisely  the existence of pre-semigeodesic charts was proved in the case when the components of the affine connection are real analytic functions.
In the proof  S.~Kowa\-lew\-sky's Theorem \cite{kowa} was used.\smallskip

We proved that the pre-semigeodesic charts exist in the case when the components of the affine connection are twice differentiable functions.  The  following is true
\begin{theorem}\label{th3}
For any affine connection determined by $\G^h_{ij}(x) \in C^r(U)$, $r \geq 2$, there exists a local transformation
of coordinates determined by $x'= f(x) \in C^{r}$ such that the connection in the new coordinates $(U',x')$, $U'\subset U$, satisfies 
$\G'{}^h_{11}(x')=0$, $h = 1,\dots, n$, 
i.e. the coordinates $(U',x')$ are pre-semigeodesic and the
components \mbox{$\G'{}^h_{ij}(x')\in C^{r-2}(U')$}.
\end{theorem}

The differentiability class $r$ is equal to $0, 1, 2, \dots ,\infty,\o$, where 0, $\infty$ and~$\o$ denotes continuous, infinitely differentiable, and real analytic functions, respectively.

It therefore follows that the existence of a pre-semigeodesic chart is guaranteed in the case when the components of the affine connection $\nabla$ are twice differentiable. 
The existence of this chart is not excluded in the case when the components are only continuous.

The components \mbox{$\G'{}^h_{ij}(x')$} can have better differentiability than $C^{r-2}(U')$. 
On the other hand, if the transformation  $x'= f(x) \in C^{r^*}$, $2\leq r^*\leq r$, leads to pre-semigeodesic coordinates (which  is possible), then it guarantees that \mbox{$\G'{}^h_{ij}(x')\in C^{r^*-2}$}.\medskip

The affine connection $\nabla$ is defined in general coordinates by 
$n^3$ components $\Gamma^h_{ij}(x)$ which are functions of $n$ variables, and $\nabla$ without torsion is defined by $n^2(n+1)/2$ components.

Theorem 3 implies that in pre-semigeodesic coordinates the number of independent functions, which are defined by $\nabla$, is reduced by $n$ functions.
It follows that all affine connections $\nabla$ in dimension $n$ depend locally only on $n(n^2- 1)$ arbitrary functions of $n$ variables, and all affine connections without torsion depend only on $n(n-1)^2/2$ arbitrary functions of $n$ variables.
  
	A manifold \an\ with a symmetric affine connection is called an 
	\emph{equiaffine manifold} if the Ricci tensor is symmetric, or equivalently, in any local coordinates $x$ there exists a function $f(x)$ satisfying 
	\cite{metal,mvh,si,no}: $$\G^\a_{i\a}=\partial_if(x).$$	
	It is clear to see that for equiaffine connections  the number of these functions is reduced by further $(n-1)$ functions.

\section{Special coordinates generated by vector fields}
\label{sec2}

Let $X$ be a vector field which is defined in the neighbourhood of the point $p$ on an $n$-dimensional manifold~$M_n$ in the coordinate system 
$x=(x^1,x^2,\ldots,x^n)$ by the components $\xi^h(x)$; $\xi^h(p)\not\equiv 0$.

It is known, see \cite{ei1,ei2,mvh,metal,ne,si,pe,po,ya}, that it is possible to find a coordinate system $x'=(x'^1,\ldots, x'^n)$ such that
\begin{equation}\label{f1}
\xi'^h(x')=\delta_1^h,
\end{equation}  
where $\delta_i^h$ is the Kronecker symbol.

The coordinate transformation from $x^h$ to $x'^h$ has the form
\begin{equation}\label{f2}
x'^h=x'^h(x^1,x^2,\ldots,x^n)
\end{equation}
for which  the law of change of the components of contravariant vectors holds:
\begin{equation}\label{f3}
\xi'^h(x')=\xi^{\alpha}(x)\cdot\partial_{\alpha}x'^h(x).
\end{equation} 

This task is solved by finding solutions $ f(x)$ and $F(x)$ of the linear partial differential equations
\begin{equation}\label{f4}
\xi^{\alpha}(x)\cdot\partial_{\alpha} f(x)=0,
\end{equation}  
and
\begin{equation}\label{f5}
\xi^{\alpha}(x)\cdot\partial_{\alpha} F(x)=1,
\end{equation}  

It is known that equation \eqref{f4} has $(n-1)$ functionally independent solutions
 \begin{equation}\label{f6}
 f^2(x), f^3(x),\ldots ,  f^n(x),
\end{equation}
which are the first integrals of the system of ordinary differential equations
\begin{equation}\label{f7}
\frac{dx^h(t)}{dt}=\xi^h(x^1(t),x^2(t),\dots,x^n(t)), \quad h=1,2,\dots,n.
\end{equation}  
Equation \eqref{f5} is solved in  the same way, its solution is denoted by $ f^1(x)$.
  
	Then the searched transformation \eqref{f2} has the following form 
	\begin{equation}\label{f8}
x'^h= f^h(x).
\end{equation}

The above solution was found for $\xi^h(x)\in C^1$, 
see \cite{ei1,ei2,ne,si,pe,ya}.
\medskip

By detailed analysis, based on the \emph{Theorem of existence of the general solution and integrals} in \cite[p.~306]{erug}, the system of ordinary differential equations   \eqref{f7} has the solutions for $\xi^h(x)\in C^0$, and, the functions $\xi^i(x)/\xi^1(x)$, $i=2,\dots,n$, satisfy Lipschitz conditions.
In this case there exist $(n-1)$ functional independent integrals 
$f^i(x)\in C^1$, $i=2,\dots,n$, which are solutions of equation  \eqref{f4} in a neighbourhood of the point $p$.

Moreover, from the differential equation  \eqref{f4} we can see that for
$\xi^h(x)\in C^r$,  $i=2,\dots,n$, there exist integrals 
$f^i(x)\in C^{r}$, $i=2,\dots,n$. 
A similar statement holds for equation \eqref{f5}, i.e. $ f^1(x)\in C^{r}$.

From the above follow.
\begin{proposition}\label{pro1}
If $\xi^h(x)\in C^r$, $r\geq1$, then there exist functionally independent solutions of \eqref{f4} and \eqref{f5}: $$ f^h(x)\in C^{r}, \ \ \ h=1,2,\dots,n.$$
\end{proposition}

\begin{theorem}\label{th4}
Let  $X$ be a vector field on $M_n$ such that
$X_p\neq0$ at a point $p\in M$. If $\xi^h(x)\in C^r,\, r\geq 1$, then there is a coordinate system $x'$ near $p$ such that $X = \partial/\partial x'^1$ and the transformation $x'=f(x)\in C^{r}$.
\end{theorem}

\begin{remark}\rm
The proof for $X{\in}C^1$ can be given e.g. by means of local flows \cite{ne}.
\end{remark}

\begin{remark}\rm
It is easy to show examples where $\xi^h(x)\in C^r$, $r\geq1$, and solutions $f^i(x)\in C^{r+1}$ do not exist.
\end{remark}

\begin{remark}\rm
Finally, we show another approach to the transformation $x'{}^h=f^h(x)$ of  Theorem \ref{th4}.

Let $\xi'{}^h(x')\in C^r$ be a vector field in coordinates $x'$ and $x'{}^h=x'{}^h(x)$ be the transformation of coordinates $x\mapsto x'$ for which $0\mapsto 0$.
Further we assume that $\xi^h(x)=\delta^h_1$. Then formula \eqref{f3} has the following form
\begin{equation}\label{f9}
\partial_1x'{}^h(x)=\xi'{}^h(x'(x)).
\end{equation}


We can look at the partial differential  equations \eqref{f9} as ordinary differential equations in the value $x^1$ and real parameters $\tilde x=(x^2,\dots,x^n)$:
\begin{equation}\label{f10}
dx'{}^h(x^1,\tilde x)/dx^1=\xi'{}^h(x'(x^1,\tilde x)),
\end{equation}
and we can use the integral form:
\begin{equation}\label{f11}
x'{}^h(x^1,\tilde x)=\varphi^h(\tilde x)+
\int^{x^1}_0
\xi'{}^h(x'(\tau^1,\tilde x))\,d\tau^1,
\end{equation}
where $\varphi^h(\tilde x)$ are functions. These functions are initial conditions for the differential equations~\eqref{f10}. For these conditions we assume
\begin{equation}\label{f12}
x'{}^h(0,\tilde x)=\varphi^h(\tilde x), \ \  h=1,2,\dots,n.
\end{equation}
Evidently, the points $(x^1,\tilde x)$ belong to a certain neighbourhood of the origin 0.

As it is known \cite{erug,ka2}, if $\varphi^h$ and $\xi'{}^h$ are continuous, then equations \eqref{f11} (and also equation \eqref{f10} with initial conditions \eqref{f12})
has a solution $x'{}^h(x^1,\tilde x)$.
For this solution, evidently, exists the partial derivative $\partial_1x'{}^h(x)$, unfortunately in general $\partial_ix'{}^h(x)$, $i=2,\dots,n$, can not exist, and in this case $x'{}^h(x)\not\in C^1$.

From properties of integrals and convergence of series of functions with parameters, see 
\cite[p.~300]{kudr}, after differentiation of the integral equations \eqref{f11} we obtain that
$$
\hbox{if \ \ }\xi'{}^h(x'), \varphi^h(\tilde x)\in C^r \ (r\geq0)
\hbox{\ \ \ then \ \ } x'{}^h(x)\in C^r. 
$$
Often it can happen that $x'{}^h(x)\in C^{r+1}$.

We note that for the transformation of coordinates  $x'(x)$ the  initial functions $\varphi(\tilde x)$ must satisfy the following conditions
$\det\|\partial_ix'{}^h(0,\tilde x)\|\neq0$ where $(0,\tilde x)$ is in neighbourhood of the origin 0.  These conditions might be for example:
$\varphi^1(\tilde x)=0$ and $\varphi^i(\tilde x)=x^i$, $i>1$.

This is the correct proof of Theorem \ref{th4}.
\end{remark}

\begin{proposition}\label{pro2}
The above general transformations $x'=f(x)$ depend on $n$~functions with $(n-1)$ arguments.
\end{proposition}

\noindent{\it Proof.}
The general solution $f$ of the homogeneous equation \eqref{f4} is a functional composition of $(n-1)$ independent solutions of~\eqref{f6}:
$f= \Phi( f^2, f^3,\ldots, f^n)$.
 The same holds for the solution of equation \eqref{f5}, because a general solution of the non-homogeneous equation \eqref{f5} is a sum of one solution of \eqref{f5} and the general solution $f$ of the homogeneous equation \eqref{f4}.
\smallskip

From the above follows that the functions $f^i$ of the transformation \eqref{f8} can have a lower class of  differentiability than  $C^{m+1}$,
it depends on the differentiability of the functions $\Phi$.

In the law of the transformation the components of the transformed tensor  depend of the  components of the tensor $T$ and also on $\partial_ix'^h$ (or $\partial'_ix^h$).
From that follows that  the introduced coordinate transformation 
$f$: $x\mapsto x'$ belongs to the class of differentiability $C^{r+1}$ 
the components of the tensor fields $T(x)\in C^{r^*}$ are transformed by 
$$T^{\cdots}_{\cdots}(x)\in C^{r^*} \longmapsto T'^{\cdots}_{\cdots}(x')\in C^{\min\{r^*,r\}}. $$
Because the transformation  law of the affine connection  \eqref{afflaw} contains  $\partial_{ij}x'^h$
$$\Gamma^h_{ij}(x)\in C^{r^*} \longmapsto \Gamma'^h_{ij}(x')\in C^{\min\{r^*,r-1\}}. $$

\section{Pre-semigeodesic coordinates}

Let $A_n=(M,\nabla)$ be an $n$-dimensional manifold $M$ with affine connection $\nabla$, and let $U\subset M$ be a coordinate neighborhood at the point $x_0\in U$. $(U,x)$ are coordinate maps on $A_n$.

It is well known that the curve $\ell$\,: $x^h=x^h(\tau)$ is a {\it geodesic}, if on it exists a parallel tangent vector. 
A geodesic $\ell$ is  characterized by the following equation 
$\nabla_{\lambda(\tau)}\lambda(\tau)=0$, where $\lambda(\tau)=dx^h(\tau)/d\tau$, ($\tau$ is  a natural parameter on $\ell$), which we can rewrite in  local coordinates 
\begin{equation}\label{geo-t}
\frac{d^2x^h(\tau)}{d\tau^2}+\Gamma^h_{ij}(x(\tau))\ 
\frac{dx^i(\tau)}{d\tau}\ \frac{dx^j(\tau)}{d\tau}
=0.
\end{equation}

The coordinates in \an\ are called
\emph{pre-semigeodesic  coordinates} if one system  of coordinate lines are geodesics and their natural parameter is just the first coordinate, see Definition \ref{def1}.

Let the $x^1$-curves be  geodesics $\ell=(\tau,x^2_0,x^3_0,\ttt,x^n_0)$ where
$\tau$ is a natural parameter.
Substituting this parametrization into the equations for geodesics \eqref{geo-t} we obtain
\be{pre-semi} \G^h_{11}(x)=0. \ee

This condition is necessary and sufficient  for a coordinate system to be pre-semigeodesic, see
 Theorem \ref{th2}.
\smallskip
\noindent{\it Proof} of Theorem \ref{th2}.
Let $\G^h_{11}=0$ hold  for $h=1,\dots ,n$. Then the local curves with parametrizations 
$\ell=(\tau,x^2_0,x^3_0,\ttt,x^n_0)$ satisfy
\be{pre-semi3}
d\ell(\tau)/d\tau=(\partial_1)_{\ell(\tau)}, \qquad d^2\ell(\tau)/d\tau^2=0,
\ee
therefore they are solutions to the system \eqref{geo-t}.

Conversely, if the curves $\ell=(\tau,x^2_0,x^3_0,\ttt,x^n_0)$ are among the solutions to \eqref{geo-t}, then due to \eqref{pre-semi3}, we get $\G^h_{11}=0$.

Hence the  pre-semigeodesic chart is fully characterized by the condition~\eqref{pre-semi} that the curves $x^1=\tau$, $x^i=\mbox{\rm const}$,
$i=2,\dots, n,$ are geodesics of the given connection in the coordinate neighbourhood, see Theorem \ref{th1}.
The definition domain $U$ of such a chart is ``tubular", a tube along geodesics.

\section{On the existence of pre-geodesic charts}

We proved Theorem \ref{th3} that a pre-semigeodesic chart exists in the case if the components of the connection are twice differentiable. 

Evidently, the existence of this chart is not excluded in the case when the components are  only continuous.

\bigskip

Let $(U,x)$ be a coordinate system at a point $p\in U\subset M$, and let
\mbox{$\G^h_{ij}(x)\in C^r$}, $r\ge0$, be components of $\nabla$ on $(U,x)$.

In a neighbourhood of $p$ we construct a set of geodesics, which go through the point $x_0=(x^1_0, x^2_0,\ldots ,x^n_0)$ of a hypersurface $\sigma\ni p$ in the direction $\lambda_0(x_0)\neq0$,  which is not tangent to $\sigma$.

Let $\sigma$ and $\lambda_0$ be defined in the following way:
\be{xx3.1}
\sigma\,{:}\ x^1=\varphi(x^2_0,\ttt,x^n_0),\ x^i=x^i_0,\ i>1,
\hbox{ \ \ and \ \ }
\l_0^h=\Lambda^h(x^2_0,\ttt,x^n_0).
\end{equation}

Then the above considered geodesics are the solutions of the following ODE's 
\be{xx3.2}
\begin{array}{l}\ds
\frac{d x^h(\tau)}{d \tau}= \l^h(\tau),\\[4mm]\ds
\frac{d \l^h(\tau)}{d \tau}=-
\G^h_{\a\b}(x(\tau))\l^\a(\tau) \l^\b(\tau)
\end{array}
\end{equation}
for the initial conditions
\be{xx3.3}
\begin{array}{ccl}
x^h(0)&=&(\varphi(x^2_0,\dots,x^n_0),x^2_0,\dots,x^n_0 ),\\[3mm]
\l^h(0)&=&\Lambda^h(x^2_0,\dots,x^n_0 )
\end{array}
\end{equation}
for any $(x^2_0,\dots,x^n_0 )$ in the neighbourhood of $p$.\smallskip

\begin {remark}\label{rem4}\rm
From \eqref{xx3.2}, \eqref{xx3.3} 
and from the theory of ODE's \cite{ka2,korn}  follows:\\
1) If $\G^h_{\a\b}(x)$ are continuous, then by the Peano existence theorem locally exists a solution. \smallskip
2) If $\G^h_{\a\b}(x)$ satisfy Lipschitz conditions, then by the 
Picard-Lindel\"of theorem this solution is unique. \smallskip
\end{remark}
In the neighbourhood of $p$ we have constructed a vector field $\l^h(x)\neq0$ which is tangent to the considered geodesics.

In addition, by more detailed analysis it can be shown that $\l^h(x)\in C^r$ if 
$\G^h_{ij}(x)\in C^r$ and moreover 
$$
\varphi(x^2,\ttt,x^n)\in C^r \hbox{ \ \ and \ \ }\Lambda^h(x^2,\ttt,x^n)\in C^r.
$$
Note 
that from the decreasing of the degree of differentiability of the functions $\varphi$ and $\Lambda^h$  follows the  decreasing of the degree of differentiability of $\l^h(x)$.

As an example, we can take the initial conditions \eqref{xx3.3} in the form:
\be{xx3.4}
x^h(0)=(0,x^2_0,\ttt,x^n_0)\hbox{ \ \ and \ \ }\l^h(0)=\delta^h_1.
\ee

Theorem \ref{th4} ensures the existence of a coordinate system $x'$ in which\linebreak 
\mbox{$\lambda'^h(x')=\delta^h_1$}. So, this system $x'$ is pre-semigeodesic, 
according to Theorems \ref{th1} and \ref{th2},
there also exists a transformation $x'=x'(x)\in C^{r}$.

 The components of a connection $\nabla$ satisfy the well-known transformation law \cite{ei2,mvh,no,metal,si}:
\begin{equation}\label{afflaw}
\G'^h_{ij}(x')=\left(
\G^\g_{\a\b}(x(x'))\ \frac{\partial x^\a}{\partial x'^i}\ 
\frac{\partial x^\b}{\partial x'^j} +
\frac{\partial^2 x^\g}{\partial x'^i\partial x'^j}
\right)\
\frac{\partial x'^h}{\partial x^\g}.
\end{equation}
Evidently,
we can prove:
 $$
\Gamma_{\a\b}^{h}(x)\in C^r,C^\infty,C^\o \ \
\longmapsto\ \
\Gamma'{}_{\a\b}^{h}(x')\in C^{r-2},C^\infty,C^\o.
$$ 

Thus  Theorem \ref{th3} was proved.
\medskip

\begin{remark} \rm
Unfortunately, the existence of a solution $\l(x)\in C^0$ (if $\G^h_{ij}\in C^0$) does not ensure the existence of a
 transformation $x'=x'(x)\in C^2$, which leads to the solution of our problem. 
In this case, the conditions for the transformations of connections are not fulfilled.
\end{remark}

\begin{remark} \rm
We show a short alternative approach of the methods for finding a transformation $x'{}^h=f^h(x)$,
$0\mapsto0$, in Theorem \ref{th3}.
\end{remark}

\noindent{\it Proof of Theorem \ref{th3}.}
From \eqref{afflaw} follows formula
$$
\frac{\partial^2 x'^h}{\partial x^i\partial x^j}=
\G^\a_{ij}(x)\, \frac{\partial x'^h}{\partial x^\a}-
\G'^h_{\a\b}(x'(x))\ 
\frac{\partial x'^\a}{\partial x^i}\frac{\partial x'^\b}{\partial x^j}.
$$
We substitute from the last formula with $i=j=1$, to the conditions $\G^h_{11}(x)=0$ and $x'=f(x)$ and we get
\begin{equation}\label{3.2}
\frac{\partial^2 f^h}{\partial x^1\partial x^1}=
-
\G'^h_{\a\b}(f(x))\ 
\frac{\partial f^\a}{\partial x^1}\frac{\partial f^\b}{\partial x^1},
\ \ \ h = 1, \dots , n. 
\end{equation}

If $x^1=t$ and if the other coordinates $\tilde x=(x^2,\dots,x^n)$ are supposed as parameters the  system \eqref{3.2} is a system of ordinary differential equations with respect to the variable~$t$.

Let the  initial condition be
\begin{equation}\label{3.3}
\begin{array}{ccc}
f^h(0,\tilde x)&=&\varphi^h_0(\tilde x)\\[2mm]
\displaystyle
\frac{\partial f^h}{\partial x^1}\ (0,\tilde x)&=&\varphi^h_1(\tilde x).
\end{array}
\end{equation}

To equations \eqref{3.2} and \eqref{3.3} Remark \ref{rem4} applies.

In addition, for the transformation $x' = f(x)$  to be regular,
it is necessary that the Jacobi matrix at a point $(0,\tilde x)$ is regular. Then it is  regular in some neighborhood of the origin 0.
An example of suitable initial conditions $(\in C^\omega)$ are 
$$
\varphi_0^1(\tilde x)=0,\ \ \varphi_0^h(\tilde x)=x^h, \ h>1, \ \
\varphi_1^h(\tilde x)=1, \ \ \varphi_1^h(\tilde x)=0.
$$

Unfortunately, from the existence of a solution does not necessary follows the existence of a
 transformation, which would lead to the solution of our problem. 

The solution $f^h(x)$ may not be generally differentiable variables
$x^2,\dots,x^n$. 
In order to realize the transformation of the components of the connection \eqref{afflaw} it is  necessary that the second derivative of $f^h(x)$ according to the variables $x^2,\dots,x^n$ exists.

It is known we can find the solution of \eqref{3.2} with initial conditions  \eqref{3.3} by a method of successive iterations \cite{ka2}:
\begin{equation}\label{3.4}
\begin{array}{rcl}
f^h_{\s+1}(x^1,\tilde x)&=&\varphi^h_0(\tilde x)+
\int\limits_0^{x_1} \lambda^h_\s(t,\tilde x)\,dt\ ,\\[3mm]
\lambda^h_{\s+1}(x^1,\tilde x)
&=&
\varphi^h_1(\tilde x) +
\int\limits_0^{x_1}\G'^h_{\a\b}(f^i_\s(t,\tilde x))\ 
\lambda^\a_\s(t,\tilde x)\lambda^\b_\s(t,\tilde x)\,dt.
\end{array}
\end{equation}

In the neighbourhood of the point $(0,x^2,\dots,x^n)$ the iterations
$f^h_{\s+1}(x^1,\tilde x)$ and $\lambda^h_{\s+1}(x^1,\tilde x)$ uniformly converge to the solutions $f^h(x)$ and $\lambda^h(x)$.

From the properties of the derivative of the integral of the parametric functions, see \cite[p.~300]{kudr},  it follows that the first derivative of solution $f^h(x)$ exists, if $\G'^h_{ij}(x'),\varphi^h_0(\tilde x),\varphi^h_1(\tilde x)\in C^1$. If we take $f^h_0=\lambda^h_0=0$, then 
each successive iteration $f^h_\s,\lambda^h_\s$ will belong 
to the class $C^1$. Because iteration is uniformly convergent,  and based on the above properties, the limits $f^h_\s\mapsto f^h$ and $\lambda^h_\s\mapsto\lambda^h$ also belong to class~$C^1$.

Analogically, the solution $f^h(x)\in C^r$ exists, if $\G'^h_{ij}(x'),\varphi^h_0(\tilde x),\varphi^h_1(\tilde x)\in C^r$.

\bibliography{\jobname}
\end{document}